\documentclass[12pt]{article}
\usepackage{amsmath,amssymb,amsthm,amsfonts,graphicx,latexsym,pdfsync}
\usepackage{mathrsfs,color}
\usepackage[latin1]{inputenc}
\numberwithin{equation}{section}

\DeclareMathOperator{\dive}{div}

\DeclareMathOperator{\supp}{supp}

\newcommand{\R}{\mathbb R}
\newcommand{\N}{\mathbb N}

\newcommand{\be}{\begin{equation}}
\newcommand{\ee}{\end{equation}}
\newcommand{\ds}{\displaystyle}

\def\B{{\mathcal B}}

\def\D{{\mathcal D}}
\def\F{{\mathcal F}}
\def\J{{\mathcal J}}

\def\LL{{\mathcal L}}
\def\M{{\mathcal M}}
\def\P{{\mathcal P}}

\def\W{{\mathcal W}}

\theoremstyle{plain}
\newtheorem{thm}{Theorem}[section]

\newtheorem{defi}[thm]{Definition}
\theoremstyle{remark}
\newtheorem{rem}[thm]{Remark}
\newtheorem{exmp}[thm]{Example}

\newtheorem*{ack}{Acknowledgements}
\newcommand{\bib}[4]{\bibitem{#1}{\sc#2: }{\it#3. }{#4.}}

\title{Evolution models for mass transportation problems}

\author{Giuseppe Buttazzo\textsuperscript{1}}
\date{}

\begin{document}

\maketitle

\footnotetext[1]{\scriptsize\ Dipartimento di Matematica, Universit\`a di Pisa, Largo B. Pontecorvo 5 - 56127 Pisa, ITALY \texttt{buttazzo@dm.unipi.it}}
\bigskip

\begin{abstract}
We present a survey on several mass transportation problems, in which a given mass dynamically moves from an initial configuration to a final one. The approach we consider is the one introduced by Benamou and Brenier in \cite{bebr00}, where a suitable cost functional $F(\rho,v)$, depending on the density $\rho$ and on the velocity $v$ (which fulfill the continuity equation), has to be minimized. Acting on the functional $F$ various forms of mass transportation problems can be modeled, as for instance those presenting congestion effects, occurring in traffic simulations and in crowd motions, or concentration effects, which give rise to branched structures.
\end{abstract}

\medskip\noindent
{\bf2010 Mathematics Subject Classification:} 49J45, 49Q20, 35F20, 76A99

\bigskip\noindent
{\bf Keywords:} Optimal transport, continuity equation, congested transport, branched transport, functionals on measures.

\section{Introduction}\label{sint}

Mass transportation theory goes back to Gaspard Monge: in 1781 he proposed in \cite{mo81} a model to describe the work necessary to move a mass distribution $\mu_1=\rho_1\,dx$ into a final destination $\mu_2=\rho_2\,dx$, given the unitary transportation cost function $c(x,y)$ which measures the work to move a unit mass from $x$ to $y$. The goal is to find a so-called {\it optimal transportation map} $T$ which moves $\mu_1$ into $\mu_2$, i.e. such that
$$\mu_2(E)=\mu_1\big(T^{-1}(E)\big)\qquad\hbox{for every measurable set }E,$$
with minimal total transportation cost
$$\int c\big(x,T(x)\big)\,d\mu_1.$$
The measures $\mu_1$ and $\mu_2$ have equal mass (normalized to one for simplicity) and are called {\it marginals}; the operator
$$T^\#\mu(E)=\mu\big(T^{-1}(E)\big)$$
is called {\it push forward} operator. The optimization problem then becomes
$$\min\Big\{\int c\big(x,T(x)\big)\,d\mu_1\ :\ T^\#\mu_1=\mu_2\Big\}.$$
The natural framework for this kind of problems is the one where $X$ is a metric space and $\mu_1,\mu_2$ are probabilities on $X$; however, the existence of an optimal transport map is a very delicate question, even in the classical Monge case, where $X$ is the Euclidean space $\R^d$ and $c(x,y)=|x-y|$ (see for instance \cite{su79,ampr03,evga99, cfmc02,trwa01}). Thus in 1942 Kantorovich proposed in \cite{ka42} a relaxed formulation of the Monge transport problem: the goal is now to find a probability on the product space, which minimizes the relaxed transportation cost
\be\label{mincost}
C(\mu_0,\mu_1)=\int c(x,y)\,\gamma(dx,dy)
\ee
over all admissible probabilities $\gamma$ on $X\times X$, where admissibility means that the projections $\pi^\#_1\gamma$ and $\pi^\#_2\gamma$ coincide with the marginals $\mu_1$ and $\mu_2$ respectively. The Kantorovich problem then reads
$$\min\Big\{\int c(x,y)\,\gamma(dx,dy)\ :\ \pi^\#_j\gamma=\mu_j\hbox{ for }j=1,2\Big\}.$$
The cases $c(x,y)=|x-y|^p$ with $p\ge1$ have been particularly studied, and the cost $C(\mu_0,\mu_1)$ in \eqref{mincost} provides, through the relation
$$w_p(\mu_0,\mu_1)=\big(C(\mu_0,\mu_1)\big)^{1/p},$$
the so-called {\it Wasserstein distance} $w_p$ which metrizes the weak* convergence on the space of probabilities $\P(\Omega)$. A very wide literature on the subject is available; we mention for instance the books \cite{ags08,vi03,vi09} where one can find a complete list of references.

A {\it dynamical model} of mass transportation consists in finding, given two probabilities $\rho_0$ (the initial configuration) and $\rho_1$ (the final destination), a curve $\rho_t$ of probabilities, with $t\in[0,1]$, joining $\rho_0$ to $\rho_1$ and minimizing some suitable functional. This functional should take into account the properties of the particular transport under consideration, like congestion or concentration phenomena.

Several models have been proposed to describe the various effects that may occur during a transportation path; in the next sections we will describe some of them. We want to stress the fact that in the congested situations (like for instance the crowd motion in case of panic) people move looking for paths with a low mass density, because they allow a higher velocity; we will see that this feature can be well described by a convex cost functional. On the contrary, in the cases where a mass concentration is favoured (as for instance in public transportation networks), people look for paths with a high mass density, which is related to some concavity in the cost.

\section{The path functionals approach}\label{spath}

In this section we consider a very general way to define minimizing paths for functionals defined on curves in a metric space $(X,d)$. The interesting case for dynamical models in mass transportation is when $X$ coincides with the class of all probabilities $\P(\Omega)$ endowed with the Wasserstein distance $w_p$.

Let $(X,d)$ be a metric space; for simplicity we assume that all closed bounded subsets of $X$ are compact. For every Lipschitz curve $\gamma:[0,1]\to X$ we define the {\it metric derivative} of $\gamma$ at the point $t$ as
$$|\gamma'|_X(t)=\lim_{s\to t}\frac{d(\gamma(s),\gamma(t))}{|s-t|}.$$
By Rademacher Theorem it can be seen (see \cite{amti04}) that for a Lipschitz curve $\gamma$ the metric derivative exists for a.e. $t\in[0,1]$ and that
$$d\big(\gamma(t),\gamma(s)\big)\le\int_s^t|\gamma'|_X(\tau)\,d\tau\qquad\forall s,t\in[0,1].$$

We are concerned with variational problems for functionals of the type
\be\label{funct}
\J(\gamma)=\int_0^1 J\big(\gamma(t)\big)|\gamma'|_X(t)\,dt
\ee
where $\gamma:[0,1]\to X$ varies among all Lipschitz curves with fixed endpoints $\gamma(0)=x_0$ and $\gamma(1)=x_1$.

\begin{thm}\label{exres}
Let $X$ be a metric space such that all closed bounded subsets of $X$ are compact, let $J:X\to[0,+\infty]$ be a lower semicontinuous functional on $X$, and let $x_0,x_1\in X$ be fixed. Assume also that
\begin{itemize}
\item[(H1)]the functional $\J$ is finite on at least a Lipschitz curve $\gamma_0$ joining $x_0$ to $x_1$;
\item[(H2)]the functional $J$ is bounded from below by a constant $\alpha>0$.
\end{itemize}
Then the minimum problem
\be\label{minpb}
\min\Big\{\J(\gamma)\ :\ \gamma\hbox{ Lipschitz},\ \gamma(0)=x_0,\ \gamma(1)=x_1\Big\}
\ee
admits a solution.
\end{thm}

The assumption of Theorem \ref{exres} on the metric space $(X,d)$ can be slightly relaxed and the following result holds.

\begin{thm}\label{duedist}
Let $(X,d)$ be a metric space and let $d'$ another distance on $X$ such that:
\begin{itemize}
\item[(K1)]$d'\le d$;
\item[(K2)]all $d$-bounded sets in $X$ are relatively compact with respect to $d'$;
\item[(K3)]the mapping $d:X\times X\to\R^+$ is lower semicontinuous with respect to the distance $d'\times d'$.
\end{itemize}
Let $\J$ be the functional defined in \eqref{funct}, where $J:X\to[0,+\infty]$ is assumed lower semicontinuous with respect to $d'$. Then, under the assumptions (H1) and (H2) above on $J$, for every $x_0,x_1\in X$ the minimum problem \eqref{minpb} admits a solution.
\end{thm}

\begin{rem}\label{coerci}
The coercivity assumption (H2) of Theorem \ref{exres} above can be weakened requiring only that
$$\int_0^{+\infty}\Big(\inf_{B_r(\bar x)}J\Big)\,dr=+\infty$$
for some (hence for all) $\bar x\in X$. Assumptions (H1) and (H2) above can still be weakened by simply requiring that there exists an admissible Lipschitz curve $\gamma_0$ such that
$$\J(\gamma_0)<\int_0^{+\infty}\Big(\inf_{B_r(x_0)}J\Big)\,dr.$$
\end{rem}

We refer to \cite{brbusa06} for the proofs of the results above; we apply here these results to the case when $X$ is a Wasserstein space of probabilities. More precisely, we consider a compact metric space $\Omega$ equipped with a distance function $c$ and a positive finite non-atomic Borel measure $m$. We consider the $q$-Wasserstein metric space $\W_q(\Omega)$ of all probability measures $\mu$ on $\Omega$, equipped with the $p$-Wasserstein distance (with $q\ge1$)
$$w_q(\mu_1,\mu_2)=\inf\Big(\int_{\Omega\times\Omega}c(x,y)^q\,\lambda(dx,dy)\Big)^{1/q}$$
where the infimum is taken on all transport plans $\lambda$ between $\mu_1$ and $\mu_2$, that is on all probability measures $\lambda$ on $\Omega\times\Omega$ whose marginals $\pi^\#_1\lambda$ and $\pi^\#_2\lambda$ coincide with $\mu_1$ and $\mu_2$ respectively.

We consider functions $J:\W_q(\Omega)\to[0,+\infty]$ to be used in \eqref{funct}. Functionals of this form have been studied by Bouchitt\'e and Buttazzo in a series of papers (see \cite{bobu90,bobu92,bobu93}) and it is shown that, under the weak* lower semicontinuity and a locality property, they can be represented in the integral form
$$J(\mu)=\int_\Omega f\Big(x,\frac{d\mu}{dm}\Big)\,dm
+\int_{\Omega\setminus A_\mu}f^\infty\Big(x,\frac{d\mu^s}{d|\mu^s|}\Big)\,d|\mu^s|
+\int_{A_{\mu}}g\big(x,\mu(x)\big)\,d\#(x)$$
where
\begin{itemize}
\item $m$ is a nonnegative nonatomic finite measure on $\Omega$;
\item $f:\Omega\times\R\to[0,+\infty]$ is an integrand with $f(x,\cdot)$ convex and lower semicontinuous;
\item $d\mu/dm$ is Radon-Nikodym derivative of $\mu$ with respect to $m$;
\item $\mu^s$ is the singular part of $\mu$ with respect to $m$ according to the Radon-Nikodym decomposition theorem;
\item $f^\infty$ is the {\it recession function} of $f$, defined by (the limit is independent of the choice of $s_0$ in the domain of $f(x,\cdot))$):
$$f^\infty(x,s)=\lim_{t\to+\infty}\frac{f(x,s_0+ts)}{t}\;;$$
\item $A_\mu$ is the set of atoms of $\mu$, i.e. the points $x\in\Omega$ such that $\mu(x):=\mu(\{x\})>0$;
\item $\#$ is the counting measure;
\item $g:\Omega\times\R\to[0,+\infty]$ is an integrand with $g(x,\cdot)$ subadditive and lower semicontinuous, such that $g(x,0)=0$, and fulfilling the compatibility condition
$$g_0(x,s)=\sup_{t>0}\frac{g(x,st)}{t}=f^\infty(x,s).$$
\end{itemize}

\begin{exmp}\label{examples}
Taking $f(s)=|s|^p$ with $p>1$ and $g(s)=+\infty$ (with $g(0)=0$) we have the Lebesgue type functionals
$$J(\mu)=\left\{
\begin{array}{ll}
\ds\int_\Omega|u|^p\,dm&\hbox{if }\mu=u\,dm\hbox{ with }u\in L^p(\Omega,m)\\
+\infty&\hbox{otherwise}
\end{array}\right.$$
whose domain is $L^p(\Omega,m)$. On the other hand, taking $f(s)=+\infty$ and $g(s)=|s|^\alpha$ with $\alpha<1$, we have the Dirac type functionals
$$J(\mu)=\left\{
\begin{array}{ll}
\ds\sum_k|a_k|^\alpha=\int_\Omega|\mu(x)|^\alpha\,d\#(x)&\hbox{if $\mu=\sum_k a_k\delta_{x_k}$ is a discrete measure}\\
+\infty&\hbox{otherwise}
\end{array}\right.$$
whose domain consists of discrete measures. Finally, taking $f(s)=|s|^p$ with $p>1$ and $g(s)=|s|^\alpha$ with $\alpha<1$ we have the Mumford-Shah type functionals
$$J(\mu)=\left\{
\begin{array}{ll}
\ds\int_\Omega|u|^p\,dm+\int_\Omega|\mu(x)|^\alpha\,d\#(x)&\hbox{if }\mu=u\,dm+\sum_k a_k\delta_{x_k}\\
+\infty&\hbox{otherwise}
\end{array}\right.$$
whose domain consists of measures with no Cantor part.
\end{exmp}

Putting together the characterization of weak* lower semicontinuous functionals $J(\mu)$ above, and the abstract existence Theorem \ref{exres} we have the following result (see \cite{brbusa06} for the proof).

\begin{thm}\label{stima F>c}
Suppose that $f(s)>0$ for $s>0$ and that $g(1)>0$. Then we have $J\ge c>0$, so that the coercivity condition (H2) of Theorem \ref{exres} is fulfilled. Therefore, the minimum problem \eqref{minpb} admits a solution, provided that there exists at least a Lipschitz curve $\gamma:[0,1]\to\W_q(\Omega)$, with given starting and ending points, with finite cost.
\end{thm}

It is interesting to study the situations when, given two probabilities $\mu_0$ and $\mu_1$ on an Euclidean domain $\Omega\subset\R^d$, an optimal path joining them actually exits. We analyze in a more detailed way the two cases below.
\begin{itemize}
\item $f(s)=|s|^p$ with $p>1$ and $g(s)=+\infty$ (with $g(0)=0$), which can be used to describe the congestion phenomena; the functional \eqref{funct} then becomes
$$\J_p(\gamma)=\int_0^1\Big(\int_\Omega|\gamma(t)|^p\,dx\Big)|\gamma'|_{\W_q}(t)\,dt$$
where $\W_q$ is the $q$-Wasserstein space and the expression $\int_\Omega|\gamma(t)|^p\,dx$ has to be intended as $+\infty$ when the measure $\gamma(t)$ is singular with respect to the Lebesgue measure $dx$.

\item $f(s)=+\infty$ and $g(s)=|s|^\alpha$ with $\alpha<1$, which on the contrary describes the concentration phenomena. In this case the functional \eqref{funct} takes the form
$$\J_\alpha(\gamma)=\int_0^1\Big(\int_\Omega|\gamma(t)|^\alpha\,d\#\Big)|\gamma'|_{\W_q}(t)\,dt$$
\end{itemize}

In the first case the domain of the functional $J$ is $L^p(\Omega)$ and the question is answered by the following result (see \cite{brbusa06}).

\begin{thm}\label{tcong}
If $\mu_0=u_0\LL^d$ and $\mu_1=u_1\LL^d$, with $u_0,u_1\in L^p(\Omega)$, then $\mu_0$ and $\mu_1$ can always be joined by a finite energy path, hence by an optimal energy path.\\
In addition, if $p<1+1/d$, then any pair of probabilities $\mu_0,\mu_1$ can be joined by a finite energy path, hence by an optimal energy path.\\
Finally, if $p\ge1+1/d$, every nonconstant Lipschitz path $\gamma(t)$ starting from a Dirac mass is such that $\J_p(\gamma)=+\infty$. Therefore, if $p\ge1+1/d$ no measure can be joined to a Dirac mass by a finite energy path.
\end{thm}

The second case, occurring when concentration phenomena are present, is somehow symmetric; in this case the domain of the functional $J$ is made only of discrete measures and the existence of optimal paths is solved by the following result (see \cite{brbusa06}).

\begin{thm}\label{tconc}
If $\mu_0$ and $\mu_1$ are finite sums of Dirac masses, i.e. $
\mu_0=\sum_Ia_i\delta_{x_i}$ and $\mu_1=\sum_Jb_j\delta_{y_j}$, then $\mu_0$ and $\mu_1$ can always be joined by a finite energy path, hence by an optimal energy path.\\
In addition, if $\alpha>1-1/d$, then any pair of probabilities $\mu_0,\mu_1$ can be joined by a finite energy path, hence by an optimal energy path.\\
Finally, if $\alpha\le1-1/d$, every nonconstant Lipschitz path $\gamma(t)$ starting from the Lebesgue measure is such that $\J_\alpha(\gamma)=+\infty$. Therefore, if $\alpha\le1-1/d$ no measure can be joined to the Lebesgue measure by a finite energy path.
\end{thm}

The example below shows a simple case in which an explicit computation of the optimal concentration path can be made.

\begin{exmp}\label{expl}
In the Euclidean plane $\R^2$ consider the origin $O$, the points $A=(1,a)$ and $B=(1,-a)$ with $a>0$, and take $\mu_0=\delta_O$, $\mu_1=\frac{1}{2}\delta_A+\frac{1}{2}\delta_B$. Using the concentration transport functional of Theorem \ref{tconc} with $0\le\alpha<1$, an easy calculation shows that the optimal path $\gamma(t)$ is given by
$$\gamma(t)=\begin{cases}
\delta_{(t,0)}&\hbox{for }0\le t\le\tau_\alpha\\
\frac{1}{2}\delta_{x(t)}+\frac{1}{2}\delta_{y(t)}&\hbox{for }\tau_\alpha<t\le1
\end{cases}$$
where $\tau_\alpha=\big(1-a(4^{1-\alpha}-1)^{-1/2}\big)^+$ and
$$x(t)=\Big(t,a\frac{t-\tau_\alpha}{1-\tau_\alpha}\Big),\qquad
y(t)=\Big(t,-a\frac{t-\tau_\alpha}{1-\tau_\alpha}\Big).$$
Figure \ref{fy} shows the plot of the entire path in the case $a=1$ for some values of the parameter $a$. Notice that in this case $\tau_\alpha=0$ for $\alpha\ge1/2$, while for $\alpha=0$ we recover the Steiner problem.
\begin{figure}[ht]
\centerline{\includegraphics[height=5.0cm,width=5.0cm]{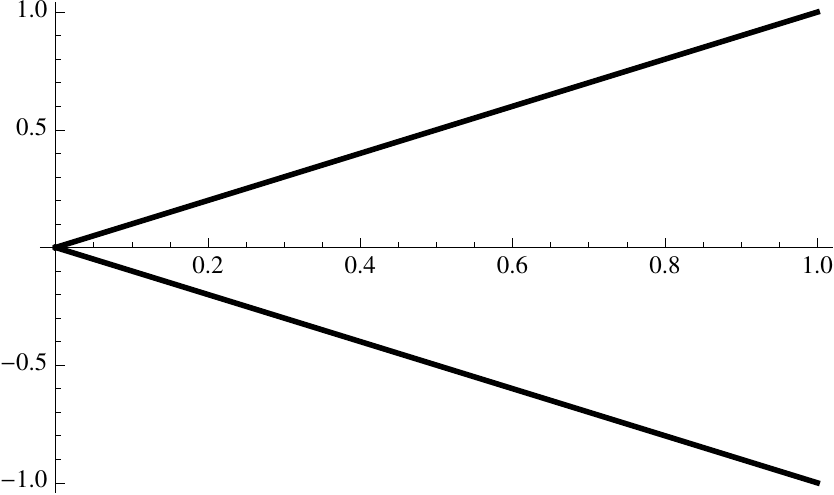}\quad\includegraphics[height=5.0cm,width=5.0cm]{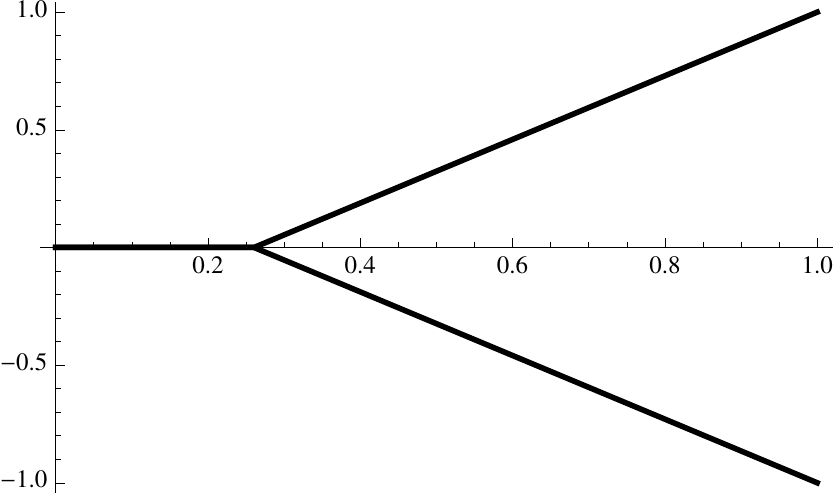}\quad\includegraphics[height=5.0cm,width=5.0cm]{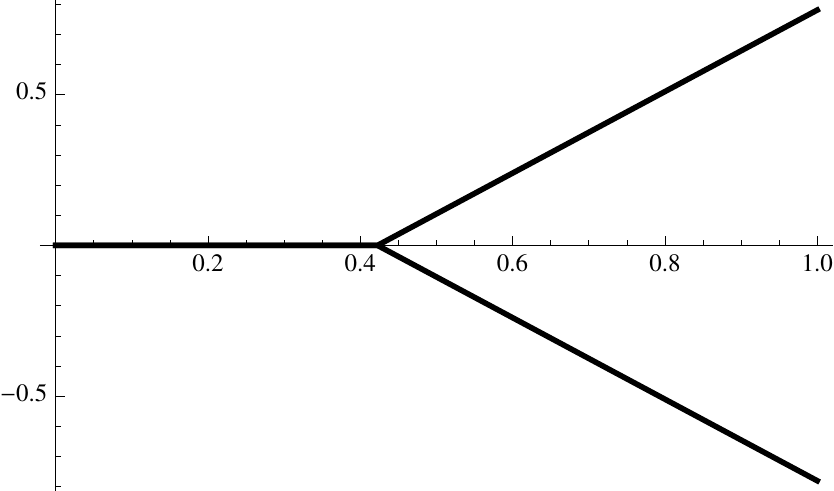}}
\caption{The optimal concentration paths for $\alpha=0.5$, $\alpha=0.25$, $\alpha=0$.}\label{fy}
\end{figure}
\end{exmp}

\section{The Eulerian approach}\label{seule}

A different dynamical approach to mass transportation problems was proposed by Benamou and Brenier in \cite{bebr00} (see also \cite{br03}). It consists in looking to pairs $(\rho_t,v_t)$ where $\rho_t$ represents the mass distribution at time $t$ and $v_t$ the related velocity field; from the mass conservation the pair $(\rho,v)$ has to satisfy the so-called {\it continuity equation}
$$\partial_t\rho+\dive_x(\rho v)=0$$
Among all pairs $(\rho,v)$ that satisfy the continuity equation above, the one providing the mass transportation, from an initial configuration $\rho_0$ to a final one $\rho_1$, is obtained by the minimization of a suitable functional $\F(\rho,v)$ that in \cite{bebr00} is taken equal to the kinetic energy:
$$\F(\rho,v)=\int_0^1\Big(\int|v|^2\,d\rho\Big)\,dt.$$
In this way the probability $\rho_0$ is dynamically transported on the probability $\rho_1$ following the {\it geodesic path} on the $2$-Wasserstein metric space $\W_2$.

To be more precise, by Theorem 8.3.1 of \cite{ags08} we have that for every absolutely continuous curve $\rho_t$ in the $p$-Wasserstein metric space $\W_p(\Omega)$ (with $p>1$) there exists a map $q$ from $[0,1]$ into the space of vector valued measures, such that $q_t\ll\rho_t$ (hence $q_t=v_t\rho_t$, being $v$ the velocity vector) which represents the flux $q=\rho v$ and satisfies
\be\label{coneq}
\partial_t\rho+\dive_xq=0\quad\mbox{and}\quad\|v_t\|_{L^p(\rho_t)}=|\rho'_t|_{W_p}\mbox{ for a.e. }t\in[0,1].
\ee
The kinetic energy is replaced in this case by the action functional
\be\label{ekin}
\F(\rho,q)=\int_0^1\Big(\int\Big|\frac{dq_t}{d\rho_t}\Big|^p\,d\rho_t\Big)\,dt.
\ee
In the degenerate case $p=1$ a little more care is needed, since the absolute continuity $q_t\ll\rho_t$ is no more guaranteed, and the $L^1$-norm has to be replaced by the mass of the measure $q_t$ (see \cite{am03}).

Note that, using the variables $\rho$ and $q$ instead of $\rho$ and $v$, provides the convexity of the functional $\F(\rho,q)$ in \eqref{ekin}.

The precise meaning of the continuity equation
$$\left\{\begin{array}{ll}
\partial_t\rho+\dive_x q=0&\mbox{ in }[0,1]\times\Omega\\
q\cdot\nu=0&\mbox{ on }[0,1]\times\partial\Omega,
\end{array}
\right.$$
has to be given in the sense of distributions, that is 
$$\int_0^1\Big[\int_\Omega\partial_t\phi(t,x)\,d\rho_t(x)+\int_\Omega D_x\phi(x,t)\cdot dq_t(x)\Big]\,dt=0$$
for every smooth function $\phi$ with $\phi(0,x)=\phi(1,x)=0$.

The general dynamical formulation of a mass transportation problems, following this Eulerian formulation, then becomes
\be\label{eeule}
\min\big\{\F(\rho,q)\ :\ \partial_t\rho+\dive_xq=0,\ \rho(0,\cdot)=\rho_0,\ \rho(1,\cdot)=\rho_1\big\}.
\ee
In the minimization problem above the continuity equation provides a linear constraint, and the existence of an optimal dynamical path $\rho_t$ easily follows by the direct methods of the calculus of variations. In the theorem below (see \cite{bujiou09}) we denote by $Q$ the time-space domain $[0,1]\times\Omega\subset\R^{1+d}$, with outer normal versor $n$, and by $\sigma$ the measures of the form $(\rho,q)$, which belong to the space $\M_b(\overline Q,\R^{1+d})$ of $\R^{1+d}$-valued measures defined on $\overline Q$. The minimization problem \eqref{eeule} can be then written in the form
\be\label{esigma}
\min\big\{\F(\sigma)\ :\ -\dive\sigma=f\hbox{ in }\overline Q,\ \sigma\cdot n=0\hbox{ on }\partial Q\big\},
\ee
where the scalar measure $f=\delta_1(t)\otimes\rho_1(x)-\delta_0(t)\otimes\rho_0(x)$ takes into account the initial-final configurations.

\begin{thm}\label{texbjo}
Let $\F:\M_b({\overline Q},\R^{1+d})\to[0,+\infty]$ be lower semicontinuous for the weak* convergence, and assume that the coercivity condition
\be\label{coercive}
\F(\sigma)\ge C|\sigma|-\frac{1}{C}\qquad\forall\sigma\in\M_b({\overline Q},\R^{1+d})
\ee
holds for a suitable constant $C>0$, where $|\sigma|$ denotes the total variation of $\sigma$ on $\overline Q$. Assume also that $\F(\sigma_0)<+\infty$ for at least one measure $\sigma_0$ satisfying the continuity equation $-\dive\sigma=f$ in $\overline Q$, with the boundary conditions $\sigma\cdot n=0$ on $\partial Q$. Then the minimum problem \eqref{esigma} admits a solution. Moreover if $\F$ is strictly convex, this solution is unique.
\end{thm}

In the case when $\F$ is convex, problem \eqref{esigma} also admits the dual formulation
\be\label{primal}
\begin{array}{lll}
(\F^*\circ A)^*(f)&=&\nonumber\min\big\{ \F(\sigma)\ :\ -\dive\sigma=f\mbox{ in }\overline Q,\ \sigma\cdot n=0\mbox{ on }\partial Q\big\}\\
&=&\ds\sup\Big\{\int\varphi\,df-\F^*(D\varphi)\ :\ \varphi\in C^1(\overline Q)\Big\},
\end{array}
\ee
where $A:C(\overline Q)\to C(\overline Q,\R^{1+d})$ denotes the operator with domain $C^1(\overline Q)$ given by
$$A(\varphi)=D\varphi\qquad\forall\varphi\in C^1(\overline Q).$$
The dual formula above holds if $\F^*$ is continuous at least at a point of the image of $A$. The primal-dual optimality condition then reads as
\be\label{epridu}
\int D\varphi_{opt}\cdot d\sigma_{opt}=\F(\sigma_{opt})+\F^*(D\varphi_{opt}),
\ee
provided an optimal solution $\varphi_{opt}$ of the dual problem exists. The point is that, in general, the maximizers $\varphi_{opt}$ of the dual problem are not in $C^1(\overline Q)$ and a relaxation procedure is necessary to make the primal-dual optimality condition meaningful. We do not deal with this rather delicate question, and we refer the interested reader to \cite{bujiou09}.

Several variants of mass transportation problems have been studied by considering in \eqref{eeule} various {\it convex} functions of the pair $(\rho,q)$.

\begin{itemize}
\item Dolbeault, Nazaret and Savar\'e considered in \cite{donasa09} functionals of the form
$$\F(\rho,q)=\int_0^1\Big(\int\Phi(\rho,q)\,dm\Big)\,dt,$$
where $m$ is a given reference measure on $\R^d$, $\rho$ and $q$ are identified through their densities with respect to $m$, and
$$\Phi(\rho,q)=\frac{|q|^p}{h(\rho)^{p-1}}=\Big(\frac{|q|}{h(\rho)}\Big)^p\!h(\rho),
\qquad p\ge1.$$
Note that, if the function $h$ is concave (for example $h(\rho)=\rho^\beta$, with $\beta\in[0,1]$), the functional $\F$ turns out to be convex as well. Functionals of this kind are motivated to provide efficient models to study diffusion PDEs of the nonlinear mobility type
$$\partial_t\rho+\dive_x\big(h(\rho)v\big)=0,$$
which can be interpreted as gradient flows of a given functional with respect to a new family of distances generalizing the ones of Wasserstein type.

\item Buttazzo, Jimenez and Oudet considered in \cite{bujiou09} a model to describe the behaviour of a crowd under some panic effects. The dynamical model is as above, with
$$\F(\rho,q)=\int_0^1\Big(\int\frac{q^2}{\rho}+c\rho^2\,dx\Big)\,dt,$$
which is a convex functional defined on Lebesgue integrable densities $(\rho,q)$.

\item Maury, Roudneff-Chupin and Santambrogio considered in \cite{mrcs10} the problem of an efficient emergency evacuation of a crowd; the target configuration $\rho_1$ is not prescribed, being replaced by an integral functional cost which takes into account the goal of the crowd. Given a room $\Omega$ with an exit $\Gamma_{out}$ the model consists in a gradient flow evolution, in the $2$-Wasserstein metric space $\W_2(\overline\Omega)$, from an initial given density $\rho_0$. The functional governing the model is
$$F(\rho)=\left\{
\begin{array}{ll}
\displaystyle\int_\Omega d(x)\,d\rho&\hbox{if $\rho\in K$}\\
+\infty&\hbox{otherwise}
\end{array}\right.$$
where $d(x)$ denotes the distance of the point $x\in\Omega$ from $\Gamma_{out}$ and
$$K=\big\{\rho\in\P(\overline\Omega),\ \rho=\rho_\Omega(x)dx+\rho_{out},\ \rho_\Omega\le1,\ \supp(\rho_{out})\subset\Gamma_{out}\big\}.$$
\end{itemize}

The situation is more delicate if we want to provide model mass trasportation models with {\it concentration effects} through an Eulerian formulation involving the continuity equation. This because we have to consider in this case concave functionals on the space of measures, taking into account the concentration effects. This happens when the moving mass has the interest to travel together as much as possible, in order to save part of the cost; for instance this occurs in the transportation of signals along telephone cables, as discovered by Gilbert who in \cite{gil67} formulated a mathematical model for it.

The starting point is, as above, to consider the admissible class $\D$ of all pairs $(\rho,q)$ with $\rho\in C\big([0,1];\P(\Omega)\big)$ and $q\in L^1\big([0,1];\M(\Omega;\R^d)\big)$ satisfying the continuity equation
$$\left\{\begin{array}{ll}
\partial_t\rho+\dive_x q=0&\mbox{ in }[0,1]\times\Omega\\
q\cdot\nu=0&\mbox{ on }[0,1]\times\partial\Omega,
\end{array}
\right.$$
and the class $\D(\rho_0,\rho_1)$ of all pairs in $\D$ with initial and final data $\rho(0,\cdot)=\rho_0$, $\rho(1,\cdot)=\rho_1$. Our goal is to minimize on $\D(\rho_0,\rho_1)$ an integral functional of the form
$$\F(\rho,q)=\int_0^1F(\rho_t,q_t)\,dt.$$
In \cite{bbs11} we presented some natural choices for the function $F(\rho_t,q_t)$ and we showed that the integral functional $\F$ above is both lower semicontinuous and coercive with respect to a suitable convergence on $(\rho,q)$, which directly provides the existence of an optimal dynamical path.

Denoting by $G_\alpha$ ($0<\alpha<1$) the functional defined on measures
$$G_\alpha(\lambda)=\begin{cases}
\ds\int_\Omega|\lambda(\{x\})|^\alpha\,d\#(x)=\sum_{i\in\N}|\lambda_i|^\alpha&\hbox{ if }\lambda=\sum_{i\in\N}\lambda_i\delta_{x_i}\\
+\infty&\hbox{ if $\lambda$ is not purely atomic,}
\end{cases}$$
we define
$$F(\rho,q)=\begin{cases}G_\alpha(|v|^{1/\alpha}\cdot\rho)&\mbox{if }q=v\cdot\rho,\\
+\infty&\mbox{if $q$ is not absolutely continuous w.r.t. }\rho.
\end{cases}$$
In this way our functional $\F$ becomes
$$\F(\rho,q)=\int_0^1\Big[\int_\Omega|v_t(x)|\rho_t(\{x\})^\alpha\,d\#(x)\Big]\,dt
=\int_0^1\Big[\sum_{i\in\N}|v_{t,i}|\rho_{t,i}^\alpha\Big]\,dt,
\qquad(\rho,q)\in\D,$$
and the dynamical model for branched transport we consider is
\be\label{minprob}
\B_{\alpha}(\rho_0,\rho_1):=\min\big\{\F(\rho,q)\ :\ (\rho,q)\in\D(\rho_0,\rho_1)\big\}.
\ee

\begin{rem}
We point out that the weak* convergence of the pairs $(\rho,q)$ is too weak for our purposes; indeed it does not directly imply the lower semicontinuity in \eqref{minprob}, since the functional $\F$ is not jointly convex in the pair $(\rho,q)$. On the other hand, if $(\rho^n,q^n)$ is a sequence in $\D$, and we assume that
$$(\rho^n_t,q^n_t)\rightharpoonup(\rho_t,q_t),\ \mbox{ for $\LL^1$-a.e. }t\in [0,1],$$
then a simple application of Fatou's Lemma gives the desired semicontinuity property of $\F$. In fact, as a consequence of the semicontinuity of $G_\alpha$ and of the convexity of the map $(x,y)\mapsto|x|^p/y^{p-1}$, we obtain that $F$ is a lower semicontinuous functional on measures.
\end{rem}

In order to prove the existence of minimizers for problem \eqref{minprob} through the direct methods of the calculus of variations involving semicontinuity and coercivity results, we introduce a convergence which is stronger than the weak* convergence of measures on $[0,1]\times\Omega$, but weaker than weak* convergence for every fixed time $t$.

\begin{defi}
We say that a sequence $(\rho^n,q^n)$ in $\D$ $\tau$-converges to $(\rho,q)$ if
$(\rho^n,q^n)\rightharpoonup(\rho,q)$ in the weak* sense of measures and in addition
$$\sup\big\{F(\rho^n_t,q^n_t)\ :\ n\in\N,\ t\in[0,1]\big\}<+\infty.$$
\end{defi}

Note that, due to the fact that the functional $\F$ is $1$-homogeneous in the velocity, its value does not change by a reparametrizations in time. By reparametrization we mean replacing a pair $(\rho,q)$ with a new pair $(\tilde\rho,\tilde q)$ of the form $\tilde\rho_t=\rho_{\varphi(t)}$, $\tilde q_t=\varphi'(t)q_{\varphi(t)}$. This equivalently means that $\tilde q$ is the image measure of $q$ through the inverse of the map $(t,x)\mapsto (\varphi(t),x)$. The following results are obtained in \cite{bbs11}.

\begin{thm}\label{teocoe}
Let $(\rho^n,q^n)$ be a sequence in $\D$ such that $\F(\rho^n,q^n)\le C$ for a suitable constant $C$. Then up to a time reparametrization, $(\rho^n,q^n)$ is $\tau$-compact.
\end{thm}

\begin{thm}\label{teosci}
Let $(\rho^n,q^n)$ be a sequence in $\D$ $\tau$-converging to $(\rho,q)$. Then
$$\F(\rho,q)\le\liminf_{n\to\infty}\F(\rho^n,q^n).$$
\end{thm}

As a consequence we obtain the following existence result.

\begin{thm}\label{teoexi}
For every $\rho_0,\rho_1\in\P(\Omega)$, the minimization problem \eqref{minprob} admits a solution $(\rho,q)\in\D$.
\end{thm}

\begin{rem}
It has to be noticed that, similarly to what happens in the path functional approach of Section \ref{spath}, for some choices of the data $\mu_0,\mu_1$ and of the exponent $\alpha$, the functional $\F$ could be constantly $+\infty$ on every admissible path $(\rho,q)\in\D(\rho_0,\rho_1)$ joining $\rho_0$ to $\rho_1$. In \cite{bbs11} the equivalence of the Eulerian model above with other variational models for branched transportation has been proven. For these models finiteness of the minima has been widely investigated, so we can infer for instance that if $\alpha>1-1/d$ then every pair $\rho_0$ and $\rho_1$ can be joined by a path of finite energy. On the other hand, if $\alpha\le1-1/d$, $\rho_0=\delta_{x_0}$ and $\rho_1$ is absolutely continuous with respect to the Lebesgue measure $\LL^d$, then there are no finite energy paths connecting them.
\end{rem}

\begin{figure}[ht]
\centerline{\includegraphics[width=5cm]{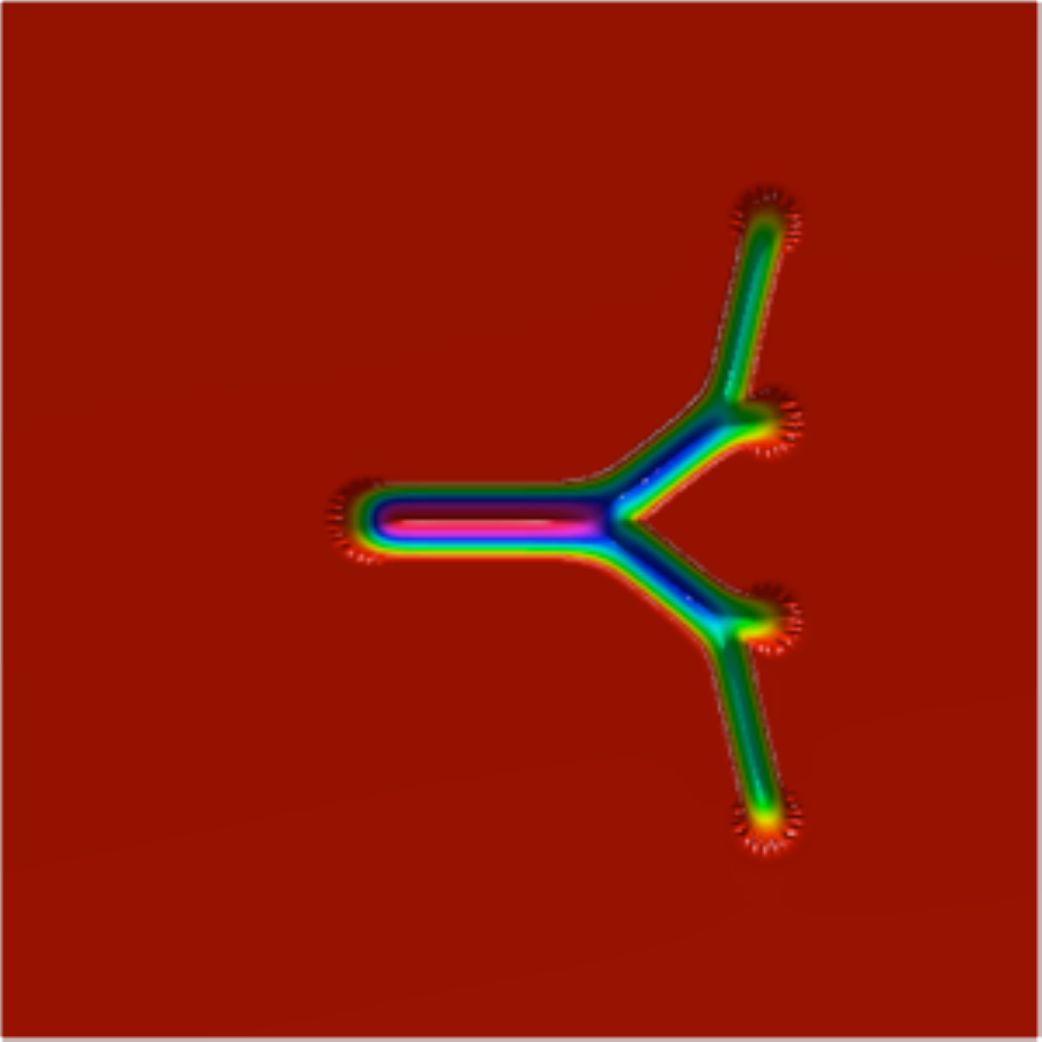}
\includegraphics[width=5cm]{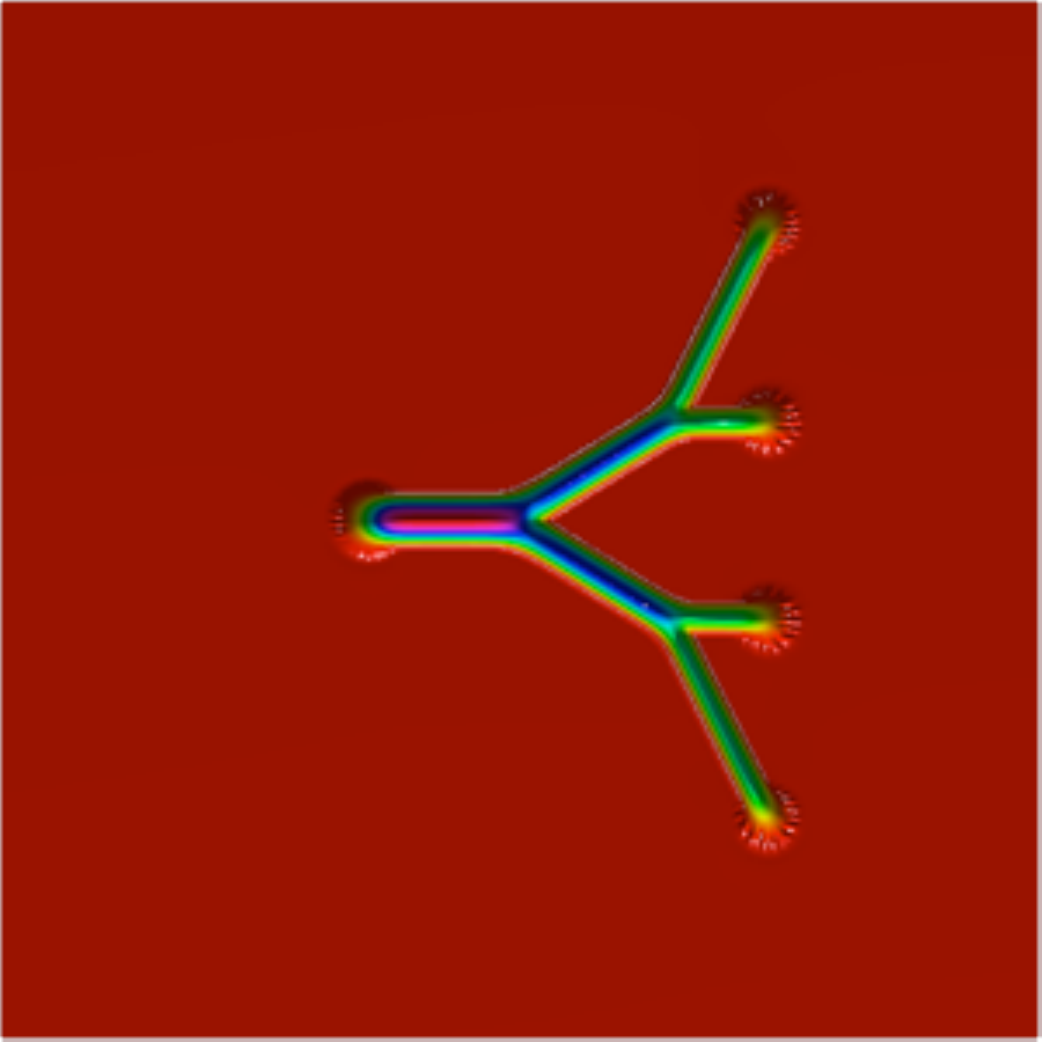}}
\centerline{\includegraphics[width=5cm]{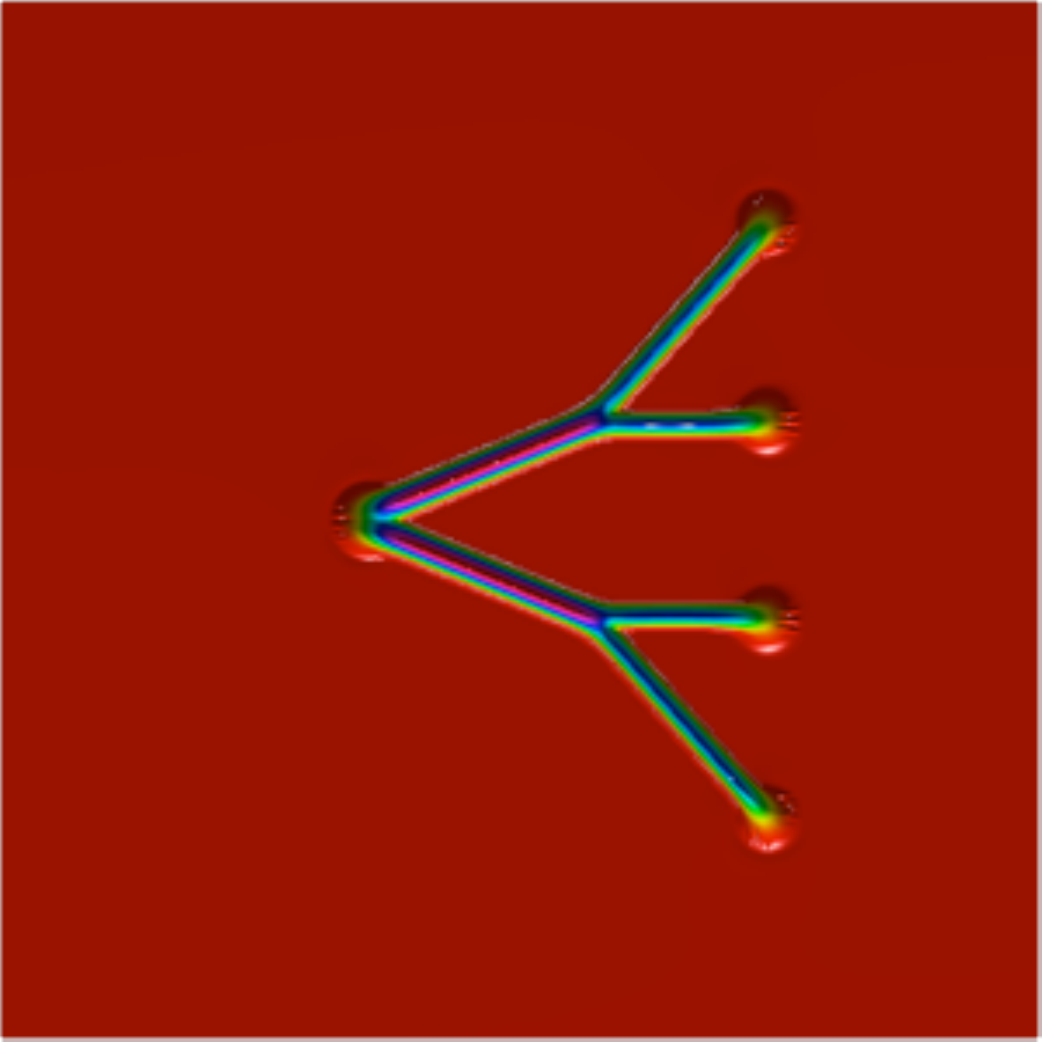}
\includegraphics[width=5cm]{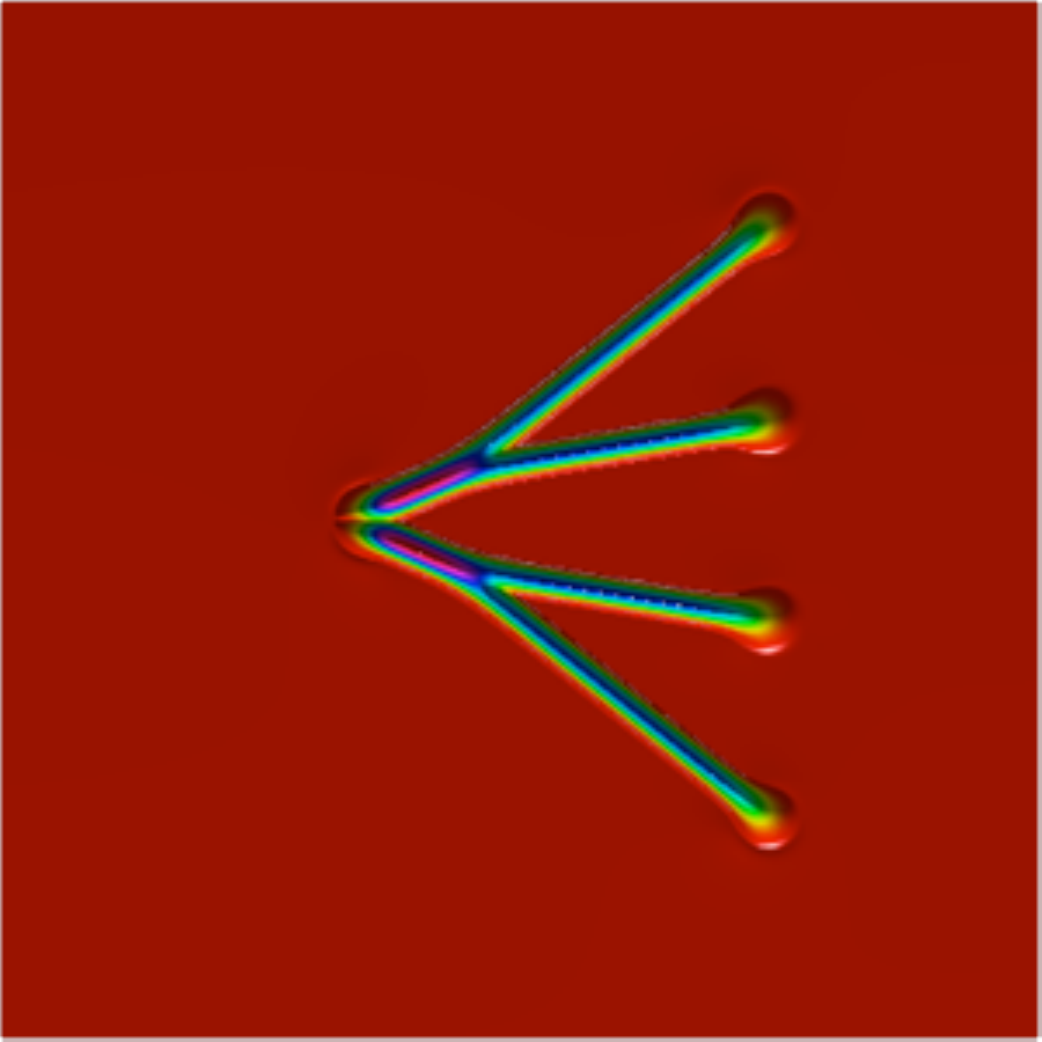}}
\caption{Optimal transport path of a point in 4 points: $\alpha=0.6,\ 0.75,\ 0.85,\ 0.95$.}
\label{foudet1}
\end{figure}

\begin{figure}[ht]
\centerline{\includegraphics[width=5cm]{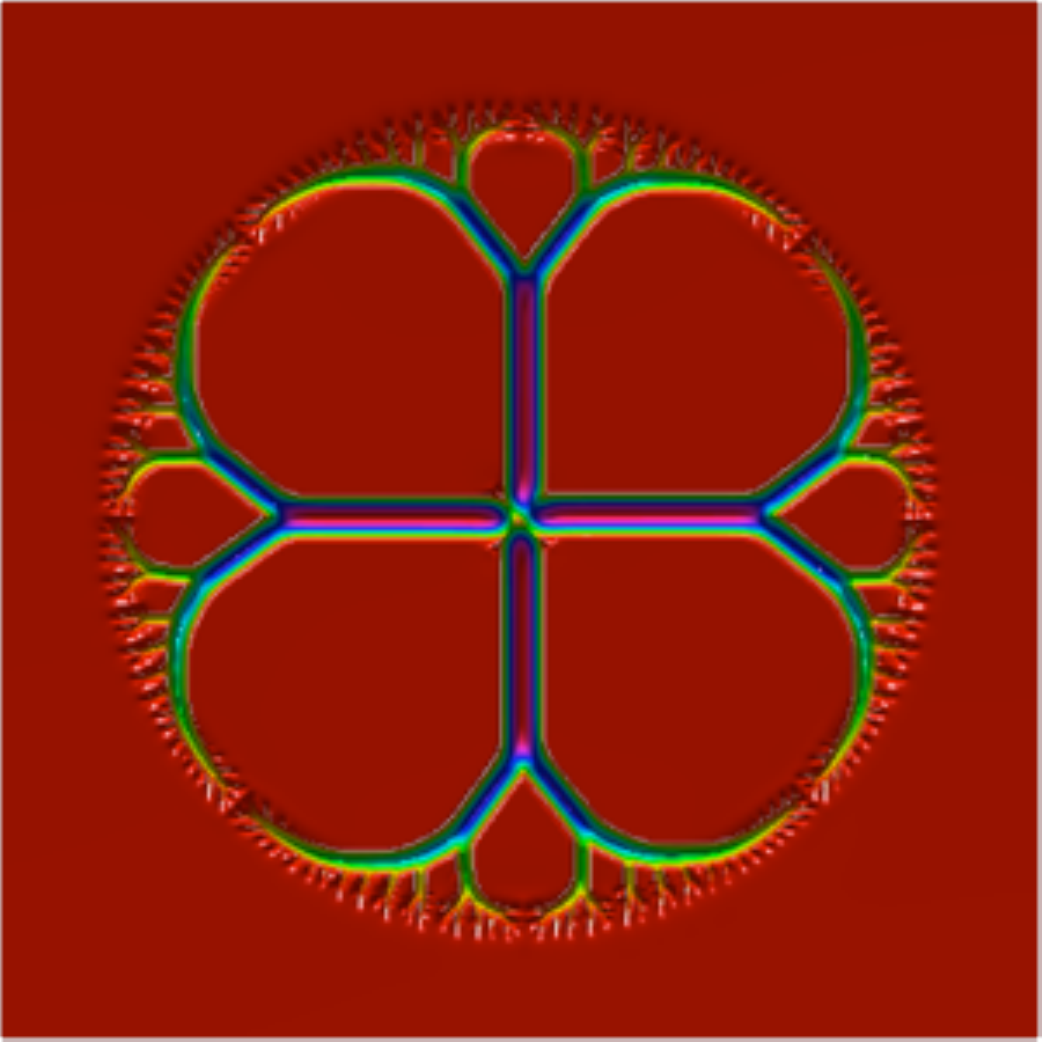}
\includegraphics[width=5cm]{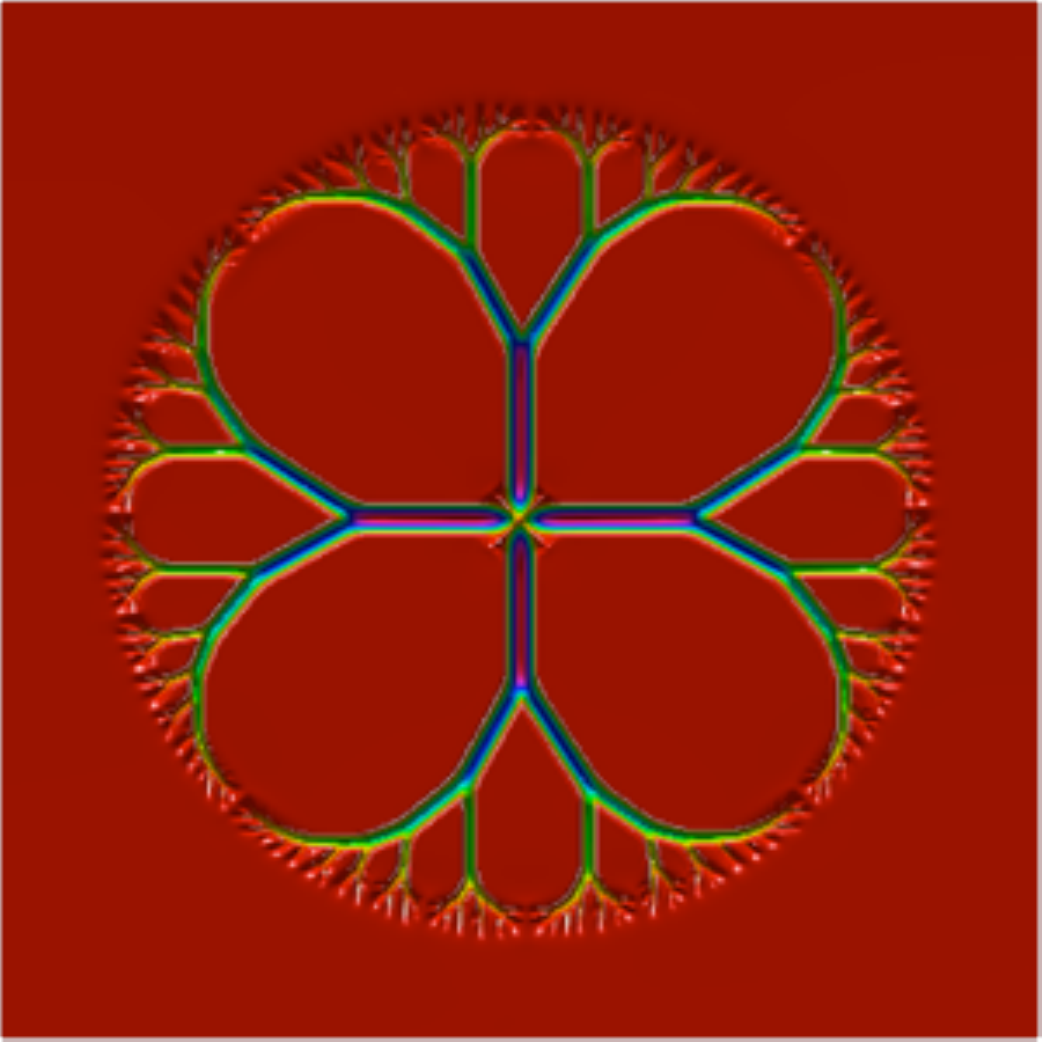}}
\centerline{\includegraphics[width=5cm]{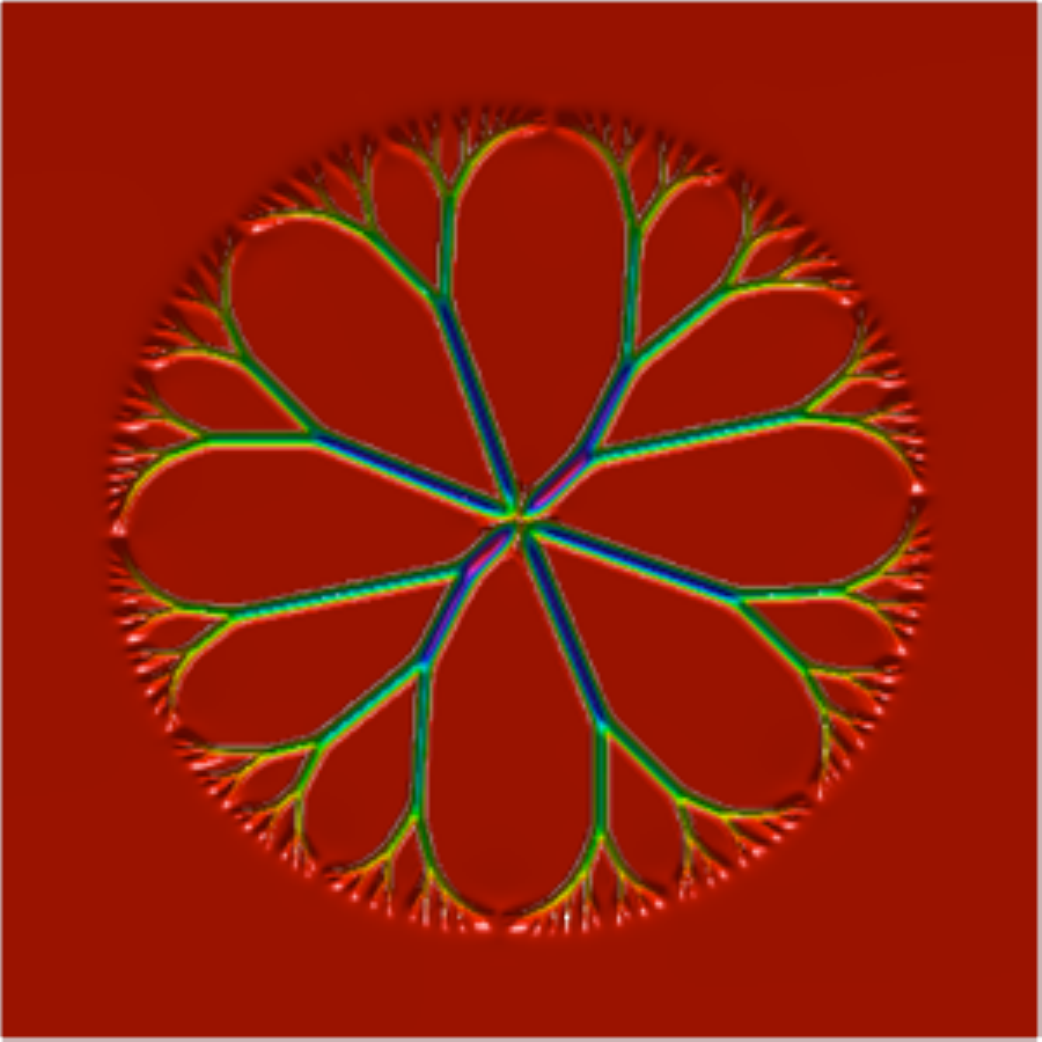}
\includegraphics[width=5cm]{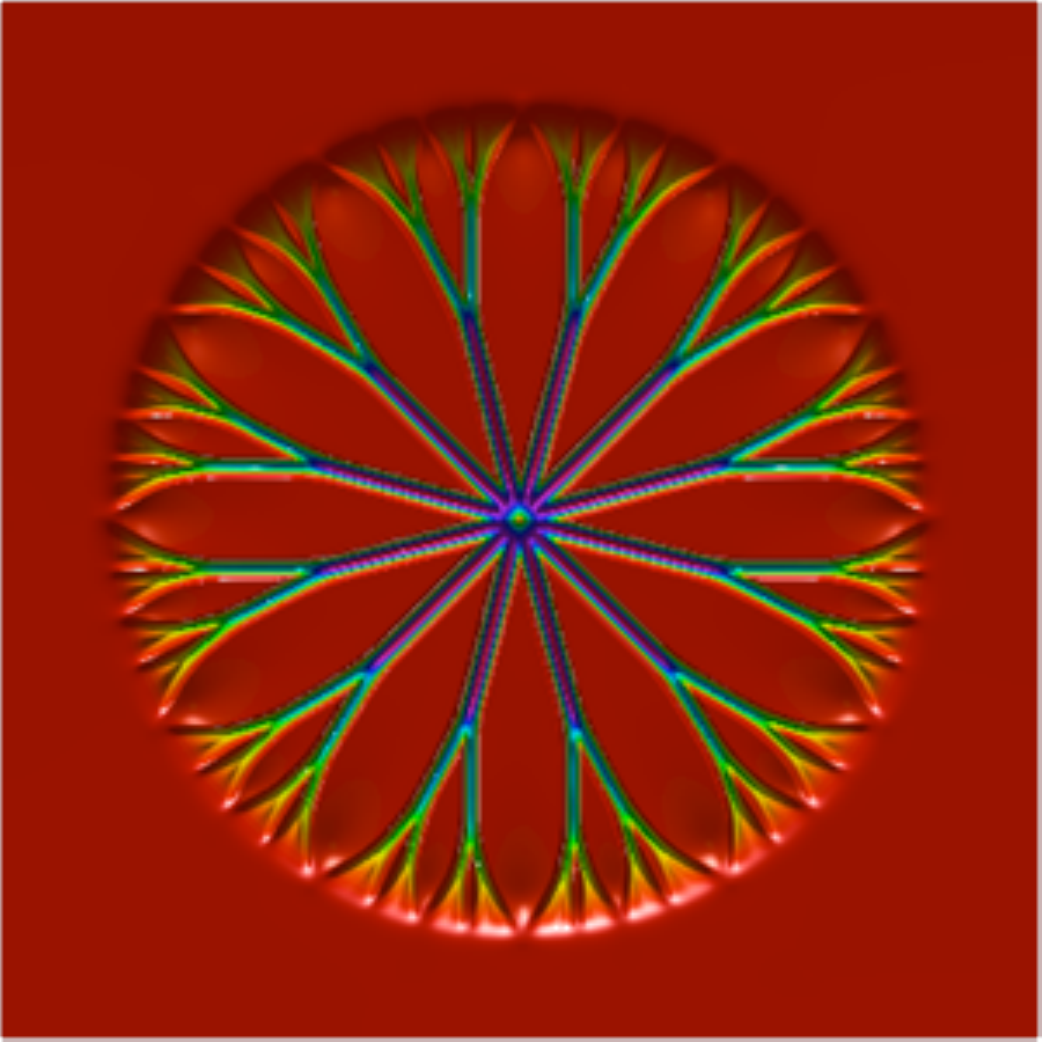}}
\caption{Optimal transport path of a point in a circle: $\alpha=0.6,\ 0.75,\ 0.85,\ 0.95$.}
\label{foudet2}
\end{figure}

Some configurations of branched transportation paths have been computed numerically by E. Oudet, providing the outputs of Figures \ref{foudet1} and \ref{foudet2} (for more examples we refer to the web site {\tt http://www.lama.univ-savoie.fr/\~{}oudet} ).

\section{Other models}\label{sother}

In addition to the evolution models presented in the previous sections other approaches are possible. We present here a model which is under investigation in \cite{buou11}. It consists in assuming that the mass is composed by many particles, each one moving under the action of a potential which takes into account the mutual interaction among them.

If $X_i(t)$ is the position of the $i$-th particle at time $t$ the Lagrangian function governing the system is
\be\label{1lagr}
L(X,X')=\frac{1}{2N}\sum_{i}|X'_i|^2+\frac{1}{N^2}\sum_{i\ne j}V(X_i,X_j)
\ee
where $N$ is the number of particles, $i,j=1,\dots,N$, and $V$ denotes the potential describing the mutual interaction among particles.

When the number $N$ of particles goes to infinity, instead of describing the motion of every single particle, it is interesting to consider the Eulerian formulation which consists in describing the movement of the mass density $\rho(t,x)$, which satisfies the {\it continuity equation}
\be\label{1conteq}
\partial_t\rho+\dive q=0.
\ee
The pair $(\rho,q)$ will then be seen as the minimizer of a suitable functional, obtained by passing to the limit as $N\to+\infty$ the Lagrangian functional \eqref{1lagr}.

The potential $V$ will of course enter in the expression of the limit functional. The term $\frac{1}{2N}\sum_{i}|X'_i|^2$ produces the kinetic energy $\int\frac{|q|^2}{2\rho}\,dx$, while the potential part $\frac{1}{N^2}\sum_{i\ne j}V(X_i,X_j)$ gives rise to the term $\int\!\!\!\int V(x,y)\rho(x)\rho(y)\,dx\,dy$. Summarizing, we end up with the energy functional
$$E(\rho,q)=\int_0^T\Big(\int\frac{|q|^2}{2\rho}\,dx
+\int\!\!\!\int V(x,y)\rho(x)\rho(y)\,dx\,dy\Big)\,dt$$
that has to be minimized among all pairs $(\rho,q)$ which satisfy the continuity equation \eqref{1conteq} together with initial and final conditions $\rho(0,\cdot)=\rho_0$ and $\rho(1,\cdot)=\rho_1$.

The energy $E(\rho,q)$ above has to be better detailed when $\rho(t,\cdot)$ is a singular measure, situation that may happen with some potential $V$. In the singular case the expression of $E(\rho,q)$ becomes
\be\label{1energy}
E(\rho,q)=\int_0^T\Big(\frac{1}{2}\int\Big|\frac{dq}{d\rho}\Big|^2\,d\rho
+\int\!\!\!\int V(x,y)\rho(dx)\rho(dy)\Big)\,dt.
\ee

Under some mild conditions on the potential $V$ it is possible to obtain an existence result for an optimal path $\rho(t,x)$. Its behaviour, in terms of congestion or concentration effects, heavily depends on the form of the potential $V$. Clearly, a potential of repulsive type will produce a congestion effect, while a potential of attractive type will produce concentration.

\begin{exmp}
In the two-dimensional case assume that the initial and final configurations are
$$\rho_0=\delta_A,\qquad\rho_1=\frac{1}{2}\delta_B+\frac{1}{2}\delta_C,$$
where $A=(0,0)$, $B=(1,1)$, $C=(1,-1)$. In the case of an attractive potential the evolution consists of the motion of two particles that travel together up to a certain point and then move symmetrically to reach their final destinations, respectively the points $B$ and $C$. If $\big(x(t),y(t)\big)$ is the path of the upper particle, we have $x(t)=t$ and $y(t)$ minimizes the functional
$$\int_0^1\big(|y'|^2+V(2y)\big)\,dt.$$
Then $y(t)=0$ for $t\le t_0$, where
$$t_0=1-\int_0^1\big(V(2y)\big)^{-1/2}\,dy,$$
while
$$y''=\nabla V(2y)\quad\hbox{for }t>t_0.$$
In the case $V(y)=C|y|^\alpha$ with $\alpha>0$ we have
$$t_0=0\quad\hbox{if }\alpha\ge2,\qquad t_0=1-\frac{2^{(2-\alpha)/2}}{(2-\alpha)\sqrt{C}}\quad\hbox{if }\alpha>2.$$
If $0<t_0<1$ we find
$$y(t)=\frac{1}{2}\big(\sqrt{C}(2-\alpha)(t-t_0)\big)^{2/(2-\alpha)}\quad\hbox{for }t>t_0.$$
Note that if $\alpha\to0$ we get
$$t_0\approx1-\frac{1}{\sqrt{C}}\qquad y(t)\approx\sqrt{C}(t-t_0)\quad\hbox{for }t>t_0$$
which shows the branching transportation behaviour.

Here are the plots of the transportation paths in some cases.
\begin{figure}[h!]
\begin{minipage}[t]{5.0cm}
\centering
\includegraphics[height=5.0cm,width=5.0cm]{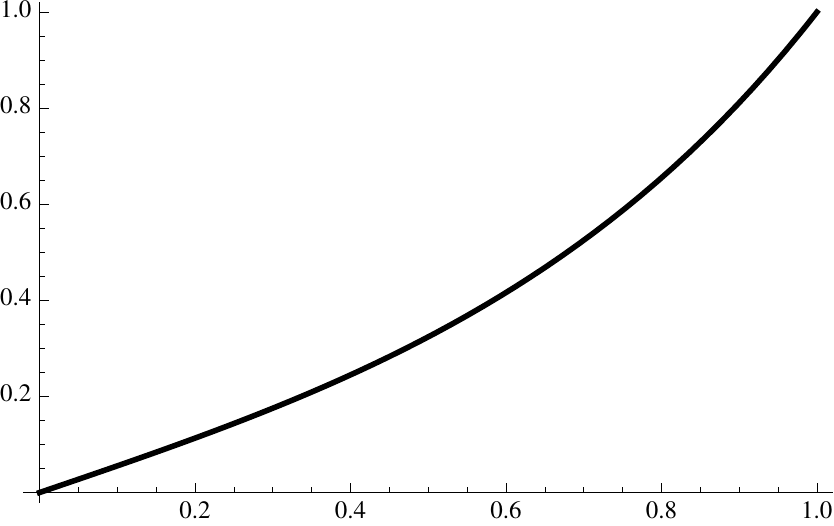}
\end{minipage}\hspace{0.1cm}
\begin{minipage}[t]{5.0cm}
\centering
\includegraphics[height=5.0cm,width=5.0cm]{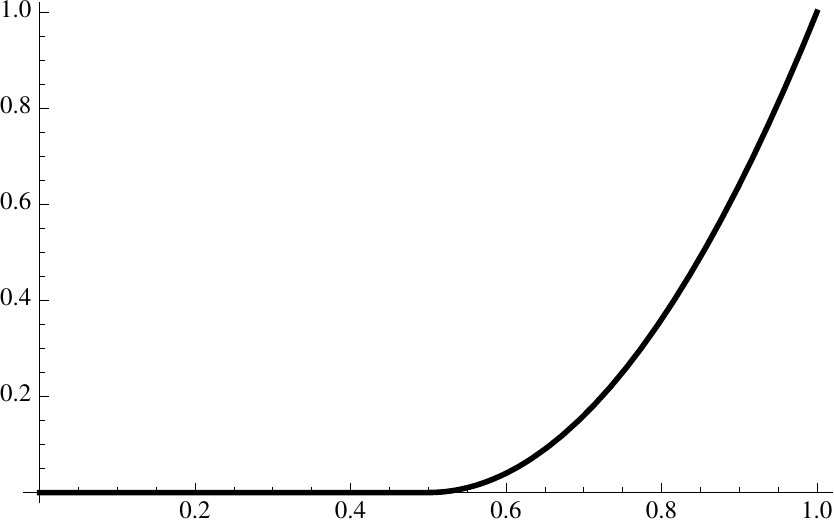}
\end{minipage}
\begin{minipage}[t]{5.0cm}
\centering
\includegraphics[height=5.0cm,width=5.0cm]{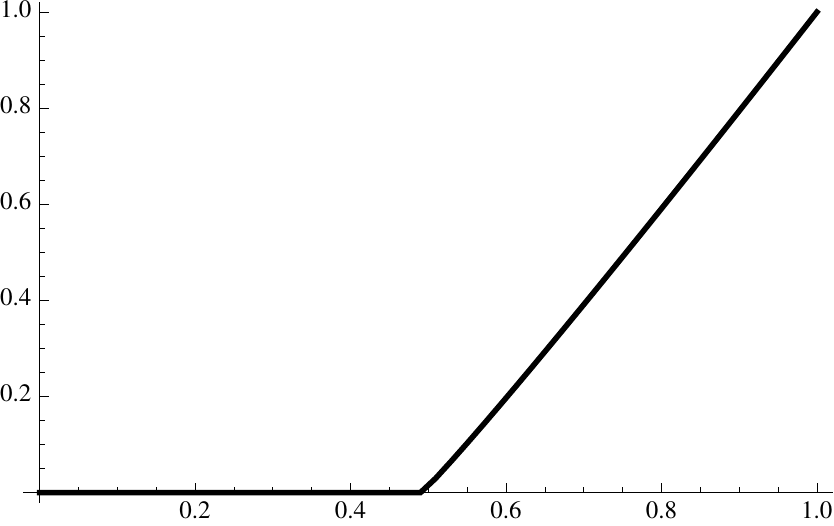}
\end{minipage}\hspace{0.1cm}
\caption{(a) $V(y)=|y|^2$.\hskip1.3truecm(b) $V(y)=8|y|$.\hskip1.3truecm(c) $V(y)=4|y|^{0.1}$.}
\label{Fig:deterministico}
\end{figure}
\end{exmp}

%It would be interesting to study, similarly to what done for the models presented in the previous sections, under which conditions on the potential $V$ every pair of probabilities $(\rho_0,\rho_1)$ can be connected by an optimal path and, when this is not possible, which kind of probabilities $(\rho_0,\rho_1)$ can be connected.

\begin{ack}
This work is part of the project 2008K7Z249 {\it Trasporto ottimo di massa, disuguaglianze geometriche e funzionali e applicazioni} financed by the Italian Ministry of Research.
\end{ack}

\end{document}